\documentclass[11pt, reqno]{amsart}

\newcommand{\mysection}[1]{\section{#1}
      \setcounter{equation}{0}}

\newcommand\supinf{\operatornamewithlimits{sup\,\,\,inf}}

\newtheorem{theorem}{Theorem}[section]
\newtheorem{lemma}[theorem]{Lemma}

\newtheorem{corollary}[theorem]{Corollary}

\theoremstyle{definition}
\newtheorem{assumption}{Assumption}[section]
\newtheorem{definition}{Definition}[section]
\theoremstyle{remark}

\newcommand\cbrk{\text{$]$\kern-.15em$]$}}
\newcommand\opar{\text{\raise.2ex\hbox{${\scriptstyle | }$}\kern-.34em$($} }
\newcommand{\tr}{\text{\rm tr}\,}

  \makeatletter
 \def\dashint{%
 \operatorname%
 {\,\,\text{\bf--}\kern-.98em\DOTSI\intop\ilimits@\!\!}}
 \makeatother

\newcommand\bR{\mathbb{R}}

\newcommand\bL{\mathbb{L}}

\newcommand\bS{\mathbb{S}}

\newcommand\bZ{\mathbb{Z}}

\newcommand\cH{\mathcal{H}}
\newcommand\cL{\mathcal{L}}

\newcommand\dist{{\rm dist}\,}

\begin{document}

\title[Rate of convergence]
{On the rate of convergence of
finite-difference approximations
for elliptic Isaacs equations in smooth domains}

\author{N.V. Krylov}
\thanks{The  author was partially supported by
 NSF Grant DMS-1160569}
\email{krylov@math.umn.edu}
\address{127 Vincent Hall, University of Minnesota,
 Minneapolis, MN, 55455}

 \keywords{Fully nonlinear
elliptic  
equations, Isaacs equations, finite differences}

\subjclass[2010]{35J60,39A14}

 \begin{abstract} 
We show that the there exists an algebraic
rate of convergence of 
solutions of 
finite-difference
  approximations for uniformly elliptic Isaacs
 in smooth bounded domains.
\end{abstract}

\maketitle

\mysection{Introduction}

Caffarelli and  Souganidis in \cite{CS} proved
that there is an algebraic rate
of convergence of solutions of finite-difference
schemes to the Lipschitz continuous viscosity solution
of the fully nonlinear elliptic equation
\begin{equation}
                                       \label{1.29.1}
F(D^{2}u(x))=f(x)
\end{equation}
in a regular domain with Dirichlet boundary data.
This is the first result available for
 fully nonlinear elliptic equations without
convexity assumptions on $F$. Naturally,
one would want to extend the result to the  $F$'s
depending also on $Du$, $u$, and $x$.
Turanova \cite{Tu_13_1} extended the results of 
\cite{CS} to $F$'s explicitly depending on $x$.
Our main goal here is to show that it is
possible by using techniques completely different from
the ones in \cite{CS} and \cite{Tu_13_1} and based on the observation
that solutions with bounded second-order
derivatives admit a H\"older estimate  of the second order
derivatives in $\cL_{d}$-sense
(see Theorem \ref{theorem 1.9.1}). However, we are only able
to deal with functions $F$ which are written in
a minimax form, that is only with Isaacs
equations, which do not encompass the most general
case of \cite{CS} and \cite{Tu_13_1}. The finite-difference approximations
we consider are also of more specific type than in \cite{CS}
and \cite{Tu_13_1}.
Our main contribution is doing away with $F$'s depending
only on the second-order derivatives.

The results of this paper bear on elliptic equations.
One can obtain similar results for parabolic Isaacs equations
(with $F$ even only measurable with respect to the time variable
using semi-discretization schemes),
because the main technical tool which is a parabolic
counterpart of Theorem \ref{theorem 1.9.1} is available
owing to Theorem 1.1 of \cite{Kr_12_2}.
Then one would ``extend'' some results in \cite{Tu_13_2}
to $F$ depending on $u,Du$, and $(t,x)$ in the same way
as the results of the present article ``extend''  the ones
in \cite{CS} and \cite{Tu_13_1}.

It is natural to compare
the results in  \cite{CS}, \cite{Tu_13_1} and in this paper
with what is known in the case
where $F$ is concave or convex with respect to $D^{2}u$.
On the one hand,  now one is able to treat
equations without convexity assumptions. But on the other hand,
in the concave or convex case  very concrete and rather good estimates
are obtained even for degenerate equations
with rates, of course, independent of the constant of
ellipticity of the equation. In  \cite{CS}, \cite{Tu_13_1} 
and in this paper
the rate is not given in any explicit way and from the proofs one
can see that this rate may depend on 
the constant of ellipticity.

The convergence of and error estimates for monotone and consistent
 approximations to
fully nonlinear, {\em first-order\/} PDEs were established a while ago by
 Crandall and Lions
\cite{CL} and Souganidis \cite{So}. 

The convergence 
of monotone and consistent
approximations for fully nonlinear, possibly degenerate 
second-order 
PDEs was first
proved in Barles and Souganidis \cite{BS}.  In a series 
of papers Kuo and 
Trudinger \cite{KT90,KT92,KT96}
  also looked in great detail at the issues of regularity and 
existence of such
approximations for uniformly elliptic equations.

There is also a probabilistic part of the story
related to controlled diffusion processes
and Bellman's equations, which
 started much earlier, see Kushner \cite{Ku},
Kushner and Dupuis \cite{KD}, also see Pragarauskas \cite{Pr}.

However, in the above cited articles
apart from \cite{CL,So}, related to the first-order equations,
and \cite{CS}, \cite{Tu_13_1}
 no rate of convergence
was established. One can read more about the past development
of the subject in Barles and Jakobsen \cite{BJ07} and  the joint article
of Hongjie Dong and the author \cite{DK07}. 
Below we concentrate  
only on some results concerning
{\em second-order\/} Bellman's equations,
which arise in many areas of mathematics such as  control
theory, differential geometry, and mathematical finance
(see Fleming and Soner \cite{FS06}, Krylov \cite{Kr77}).

The first estimates of the rate of convergence
for second-order degenerate or nondegenerate
 Bellman's equations appeared  in 1997
(see \cite{Kr97}). For equations with constant ``coefficients"
and arbitrary monotone  finite-difference approximations
it was proved in \cite{Kr97} that the rate of convergence
is $h^{1/3}$ if the error in approximating the true operators with
finite-difference ones is of order $h$ on {\em
three\/} times continuously differentiable functions. 
The order becomes
$h^{1/2}$ if the error in approximating the true operators with
finite-difference ones is of order $h^{2}$
on {\em four\/} times continuously differentiable functions
(see Remark 1.4 in \cite{Kr97}, which however contains
an arithmetical error albeit easily correctable. Also see
 Theorem 5.1 in \cite{Kr97}).
One of the   main ideas of \cite{Kr97} is that the equation and its 
finite-difference approximation should play symmetric
roles. The proofs
in \cite{Kr97} are purely analytical (in contrast with what
one can read in some papers mentioning \cite{Kr97}) 
even though sometimes probabilistic
{\em interpretation\/} of some statements are also given.
The next step was done in \cite{Kr99} where the 
so-called method of ``shaking  the coefficients''
was introduced to deal with the case of {\em degenerate\/}
parabolic Bellman's equations with {\em variable\/} coefficients.
The two sided error estimates were given: from the one side
of order $h^{1/21}$ and from the other $h^{1/3}$. Here $h$
(unlike in \cite{Kr97})
was naturally interpreted as the mesh size and
the approximating operators were assumed to approximate
the true operator with error of order $h$ on {\em three\/}
times continuously differentiable functions. 
 The order $1/21$ was improved to $1/7$
in \cite{BJ07} in the general setting of \cite{Kr99}.

These rates may look unsatisfactory and then 
one tries to get  better
estimates on the account  of using  some special approximations,
say providing the error 
  of order $h^{2}$ of approximating the main part of the
true operator
on {\em four\/} times continuously differentiable functions.
This was already mentioned in  Remark 1.4 of \cite{Kr97}
and  used by Barles and Jakobsen in \cite{BJ02}
to extend the results in \cite{Kr97} to  
equations with variable lower-order ``coefficients''.

One can also consider  special finite-difference
approximations, for instance, only containing
pure second-order differences in place of second-order
derivatives, when this $h^{2}$
approximating error is automatic. In such cases
the optimal rate $h^{1/2}$ was obtained in the joint
work of Hongjie Dong and the author \cite{DK07}
for degenerate
parabolic Bellman's equations with Lipschitz
coefficients in domains.
Both ideas of symmetry and   ``shaking the coefficients''
are used in \cite{DK07} as well as in \cite{DK05}.
In the paper by Hongjie Dong and the author \cite{DK05}
we consider among other things weakly nondegenerate
Bellman's equations with {\em constant\/} ``coefficients''
in the whole space
and obtain the rate of convergence $h$, where
$h$ is the mesh size. It may be tempting to say
that this result is an improvement of earlier
results, however it is just a better rate 
under different conditions.

 It is worth noting that
the set of equations satisfying the conditions
 in \cite{DK07} is smaller than the one in the papers by
Barles and Jakobsen \cite{BJ05,BJ07}, the results of which
obtained by using the theory of viscosity solutions
  guarantee the rate $h^{1/5}$. However, in {\em the 
concrete examples}
given  in Sections 4 in 
\cite{BJ05,BJ07}
of applications of the general scheme,  
one gets the rate $h^{1/2}$ according to \cite{DK07}
if one adds the requirement  that
the coefficients be twice differentiable
(and in the case of diagonally dominant matrices
replaces Kushner's approximation with a different one
avoiding using mixed finite differences, which is possible as shown
in \cite{Kr08}). Moreover, if the equation
is uniformly nondegenerate, then the rate is
at least $h^{2/3}$ in the elliptic case
according to the result of \cite{Kr_13_2}. 
 In \cite{BJ05,BJ07} the coefficients are only assumed to be 
once differentiable and still the rate $h^{1/5}$
is guaranteed   regardless
of degeneracy or nondegeneracy. One more point to be noted is that
in \cite{BJ07} parabolic equations are considered
with various types of approximation such as 
Crank-Nicholson and splitting-up schemes
 related to the time
derivative   combined with not necessarily
finite-difference approximation of differential
operators.

In conclusion of the introduction
we fix some notation. 
Let $\bR^{d}=\{x=(x_{1},...,x_{d})\}$
be a $d$-dimensional Euclidean space
with scalar product $\langle\cdot,\cdot\rangle$ and
let $\bS$ be the set of $d\times d$ symmetric
matrices. For a fixed constant $\delta\in(0,1]$
denote by $\bS_{\delta}$ the
subset of $\bS$ consisting
of matrices with eigenvalues in $[\delta,\delta^{-1}]$.
For a function $F(u,x)$ given
for
$$
u=(u',u''),\quad u'=(u'_{0},...,u'_{d})\in\bR^{d+1},
\quad u''\in\bS,\quad x\in\bR^{d}
$$ 
and smooth real-valued function $u(x)$ on $\bR^{d}$
we set $D_{i}u=\partial u/\partial x_{i}$, $D_{ij}u=
D_{i}D_{j}u$,
$$
F[u](x)=F(u(x),Du(x),D^{2}u(x),x),
$$
where $Du=(D_{i}u)$ is the gradient of $u$ and $D^{2}u=(D_{ij}u)$
is its Hessian.

 Let $G$ be a fixed open bounded subset 
of $\bR^{d}$
with $C^{2}$ boundary. 

\mysection{Main result}

                                             \label{section 1.30.1}

 Let $d_{1}\geq d$
  be an integer.
Assume that we are given separable metric spaces
  $A$ and $B$  and let,
for each $\alpha\in A$ and $\beta\in B$, 
  the following 
  functions on $\bR^{d}$ be given: 

(i) $d\times d_{1}$
matrix-valued $\sigma^{\alpha\beta}( x)
 =
(\sigma^{\alpha\beta}_{ij}( x))$,

(ii)
$\bR^{d}$-valued $b^{\alpha\beta}( x)=
(b^{\alpha\beta}_{i }( x))$, and

(iii)
real-valued  functions 
$c^{\alpha\beta}( x) $,   
  $f^{\alpha\beta}( x) $.

\begin{assumption}
                                    \label{assumption 1.9.1}
(i)   All the above functions are continuous with respect to
$\beta\in B$ for each $(\alpha,x)$ and continuous with respect
to $\alpha\in A$ uniformly with respect to $\beta\in B$
for each $x$. These functions are Borel measurable
functions of  $(\alpha,\beta,x)$,   and
$c^{\alpha\beta}\geq0$.

(ii) For any  $x,y\in\bR^{d}$ 
 and $(\alpha,\beta)\in A\times B $
$$
 \|\sigma^{\alpha\beta}( x )\|,|b^{\alpha\beta}( x )|,
| c^{\alpha\beta}( x)|,|  f^{\alpha\beta}( x)| \leq\delta^{-1},
$$
$$
\|\sigma^{\alpha\beta}( x)-\sigma^{\alpha\beta}( y)\|+
  |b^{\alpha\beta}(  x)-b^{\alpha\beta}( y) |\leq \delta^{-1}|x-y|,
$$
 where for a matrix $\sigma$ we denote $\|\sigma\|^{2}=
\tr\sigma\sigma^{*}$.

(iii) There exists a constant $\gamma_{1}\in(0,1]$
such that for any $x,y\in\bR^{d}$ 
 and $(\alpha,\beta)\in A\times B $   we have
$$
 |c^{\alpha\beta}( x)-c^{\alpha\beta}( y) |+
  |f^{\alpha\beta}(  x)-f^{\alpha\beta}( y) |\leq \delta^{-1}
|x-y|^{\gamma_{1}}.
$$

(iv) For $\alpha\in A$, $\beta\in B$,  
and $x \in\bR^{d}$ we have $a^{\alpha\beta}(x)\in\bS_{\delta}$,
where $a^{\alpha\beta}=(a^{\alpha\beta}_{ij})=(1/2)
\sigma^{\alpha\beta}(\sigma^{\alpha\beta})^{*}$.

\end{assumption}

We are   going to deal with
$$
H(u,x):=\supinf_{\alpha\in A\,\,\beta\in B}
[a^{\alpha\beta}_{ij}(x)u''_{ij}+b^{\alpha\beta}_{i}(x)u'_{i}
-c^{\alpha\beta}(x)u'_{0}+f^{\alpha\beta}(x)]
$$
(the summation convention over repeated indices
is enforced and here the summations are done
before other operations are performed),
which is a typical object in the theory of 
stochastic differential games.
As in \cite{Kr_14} one associates with the above objects
and zero boundary data on $\partial G$
a value function $v(x)$. According to
\cite{Kr_14} the dynamic programming principle holds, which
along with the continuity of the data allows us to use
the results of \cite{FS89} and conclude that $v$
is a unique viscosity solution of class $C(\bar G)$
of $H[v]=0$ in $G$ with zero boundary data.

Next, for smooth functions $u(x)$ and  $\alpha\in A$ and $\beta\in B$
introduce
$$
L^{\alpha\beta}u(x)=
a^{\alpha\beta}_{ij}(x)D_{ij}u(x)+b^{\alpha\beta}_{i}(x)D_{i}u(x)
-c^{\alpha\beta}(x)u(x).
$$
 
As is well known (see, for instance, \cite{KT92}),
there exists a finite set $\Lambda=\{l_{1},...,l_{d_{2}}\}
\subset\bZ^{d}$  containing all vectors from the standard 
orthonormal basis of $\bR^{d}$
such that one has the following representation
$$
L^{\alpha\beta}u(x)=
a^{\alpha\beta}_{k}(x)D_{l_{k}}^{2}u(x)+
\bar b^{\alpha\beta}_{k}(x)D_{l_{k}}u(x)
-c^{\alpha\beta}(x)u(x),
$$
where $D_{l_{k}}u(x)=\langle D u ,l_{k}\rangle$,
$a^{\alpha\beta}_{k}$ and $\bar b^{\alpha\beta}_{k}$ are certain
bounded functions and $a^{\alpha\beta}_{k}\geq\delta_{1}$,
with a constant $\delta_{1}>0$. One can even arrange for
such representation to have the coefficients 
$a^{\alpha\beta}_{k}$ and $\bar b^{\alpha\beta}_{k}$ 
with the same regularity
properties with respect to $x$ as the original ones
$a^{\alpha\beta}_{ij}$ and $b_{i}^{\alpha\beta}$
(see, for instance, Theorem 3.1 in \cite{Kr_11}).
Define
$B$ as the smallest closed ball containing $\Lambda$,
and for $h>0$ set $\bZ^{d}_{h}=h\bZ^{d}$,
$$
G_{h}=G\cap\bZ^{d}_{h},\quad
G^{o}_{h}=\{x\in \bZ^{d}_{h}:x+hB\in G\},\quad
\partial_{h}G=G_{h}\setminus G^{o}_{h}.
$$

Next, for $h>0$ we introduce
$$
\delta_{h,l_{k}}u(x)=\frac{u(x+hl_{k})-u(x)}{h} ,
$$
$$
\Delta_{h,l_{k}}u(x)=\frac{u(x+hl_{k})-
2u(x)+u(x-hl_{k})}{h^{2}} ,
$$
$$
L^{\alpha\beta}_{h}u(x)=
a^{\alpha\beta}_{k}(x)\Delta_{h,l_{k}}u(x)+
\bar b^{\alpha\beta}_{k}(x)\delta_{h,l_{k}}u(x)
-c^{\alpha\beta}(x)u(x),
$$
$$
H_{h}[u](x)=
\supinf_{\alpha\in A\,\,\beta\in B}[L^{\alpha\beta}_{h}u(x)
+f^{\alpha\beta}(x)].
$$

It is a simple fact shown, for instance, in \cite{KT92}
that   for each sufficiently small $h$ there exists a unique
  function $v_{h}$ on $G_{h}$ such that
$H_{h}[v_{h}]=0$ on $G^{o}_{h}$ and 
$v_{h}=0$ on $\partial_{h}G$.

Here is our main result:

\begin{theorem}
                                          \label{theorem 1.28.1}
There exist constants $N$ and $\beta>0$ such that
for all sufficiently small $h>0$ we have on $G_{h}$ that
$$
|v_{h}- v|\leq Nh^{\beta}.
$$

\end{theorem}

The proof of Theorem \ref{theorem 1.28.1}
is given in Section \ref{section 1.30.2}
after some preparations are made.

\mysection{A general scheme}
                                   \label{section 12.14.1}

Denote by 
$\bL_{\delta}$   the set 
of differential operators of the form
$$
Lu(x)=a_{ij}(x)D_{ij}u(x)+b_{i}(x)D_{i}u(x)+c(x)u(x),
$$
where $a(x)=(a_{ij}(x))$, $b(x)=(b_{i}(x))$,
and $c(x)$ are $\bS_{\delta}$-valued,
$\bR^{d}$-valued, and real-valued Borel
functions, respectively, on $\bR^{d}$ satisfying
$$
|b|\leq \delta^{-1},\quad 0\leq-c\leq \delta^{-1}.
$$

Here is Corollary 3.5 of \cite{Kr_13}
which slightly generalizes the main result
in Fang-Hua Lin \cite{FHL}.

\begin{lemma}
                              \label{lemma 10.7.1}
There exist  $\gamma_{0}\in(0,1)$  and $N$
depending only on $\delta$,   $d$, and $G$ such that
for any  $L\in\bL_{\delta }$,  
 $\gamma\in(0,\gamma_{0}]$, and $u\in 
W^{2}_{d,loc}(G)\cap C(\bar G)  $
 we have
\begin{equation}
                                       \label{9.5.06}
 \int_{G}(|D^{2}u|^{\gamma}+|Du|^{\gamma})
\,dx \leq N\|Lu\|_{\cL_{d}(G)}^{\gamma}
 +N\sup_{\partial G}|u|^{\gamma} .
\end{equation}

\end{lemma}

Next, we consider a function $F(u,x)$ given
for
$$
u=(u',u''),\quad u'=(u'_{0},...,u'_{d})\in\bR^{d+1},
\quad u''\in\bS,\quad x\in\bR^{d}.
$$ 

\begin{assumption}
                            \label{assumption 1.14.1}
(i) The function 
$F(u,x)$ is Lipschitz continuous with
Lipschitz constant $\delta^{-1}$
 with respect to $u$,

(ii) For any $x$ at all point of differentiability
of $F(u,x)$ with respect to $u$, we have
$F_{u''}\in\bS_{\delta}$, $F_{u'_{0}} \leq 0$,

(iii) There is a constant $\gamma_{1}
\in(0,1]$ such that for any $x,y\in\bR^{d}$ and $u=(u',u'')$
such that $|x-y|\leq1$ we have
$$
|F(u,x)-F(u,y)|\leq \delta^{-1}|x-y|^{\gamma_{1}}(1+|u|).
$$

\end{assumption}

For $\varepsilon>0$ introduce
$$
G^{\varepsilon}=\{x:\exists y\in G, |y-x|\leq 
\varepsilon\},\quad
B_{\varepsilon}=\{x:|x|<\varepsilon\}.
$$
 
The following result is one of our main technical tools.

\begin{theorem}
                                  \label{theorem 1.9.1}
Let a function $u\in W^{2}_{\infty}(G^{\varepsilon})$
satisfy
$F[u]=0$ in $G^{\varepsilon}$,
and let $p\in(\gamma_{0},\infty)$.
Then for any $h\in[0,\varepsilon\wedge 1]$ and $|l|=1$
$$
 \|u(hl+\cdot)-u\|_{W^{2}_{p}(G)} \leq
Nh^{\gamma_{1}\gamma_{0}/p} (1+M_{\varepsilon}),
$$
where 
$$
M_{\varepsilon}:=\|
u\|_{W^{2}_{\infty}(G^{\varepsilon})}
$$
and the constant
$N$ depends only on $G$, $\delta$,  $p$,
and $d$.

\end{theorem}

 Proof. Since the increments of the function and its
first derivatives are well controlled by
the product of $Nh$ times $M_{\varepsilon}$,
it suffices to prove that
\begin{equation}
                                                  \label{1.9.1}
 \|D^{2}u(hl+\cdot)-D^{2}u\|_{\cL_{p}(G)} \leq
Nh^{\gamma_{1}\gamma_{0}/p}(1+ M_{\varepsilon}).
\end{equation}
Denote $v(x)=u(hl+x)-u(x)$ and observe that
$v$ satisfies
$$
Lv+f=0
$$
in $G$,
where $L\in \bL_{\delta}$ and $f(x)=F(u(x),Du(x),D^{2}u(x),
x+hl)-F[u](x)$, so that
$$
|f|\leq \delta^{-1}h^{\gamma_{1}}(1+M_{\varepsilon}).
$$

Owing to Lemma \ref{lemma 10.7.1} for any $\lambda>0$
we have
$$
|G\cap\{|D^{2}v|^{p}\geq\lambda\}|\leq
\frac{1}{\lambda^{\gamma_{0}/p}}\int_{G}|D^{2}
v|^{\gamma_{0}}\,dx
\leq N
\frac{h^{\gamma_{1}\gamma_{0}}}
{\lambda^{\gamma_{0}/p}}(1+M^{\gamma_{0}}_{\varepsilon}).
$$
We also note that if $\lambda>(2M_{\varepsilon})^{p}$, then
$G\cap\{|D^{2}v|^{p}\geq\lambda\}=\emptyset$. It follows
that
$$
 \|D^{2}u(h+\cdot)-D^{2}u\|_{\cL_{p}(G)}^{p}
\leq N(1+M^{\gamma_{0}}_{\varepsilon})
h^{\gamma_{1}\gamma_{0}}
\int_{0}^{(2M_{\varepsilon})^{p}}
\lambda^{-\gamma_{0}/p}\,d\lambda
$$
$$
\leq
Nh^{\gamma_{1}\gamma_{0}}(1+M_{\varepsilon}^{p})
$$
and the theorem is proved.

We say that a number $\varepsilon_{0}>0$
is sufficiently small if Lemma \ref{lemma 10.7.1}
holds with the same $\gamma_{0}$,
a constant $N$ which is twice the constant
from \eqref{9.5.06}, and with $G^{\varepsilon_{0}}$
in place of $G$. The fact that the set of sufficiently
small $\varepsilon_{0}$ contains 
$[0,\alpha)$ with $\alpha>0$ follows from the way
Corollary 3.5 of \cite{Kr_13} is proved and from the fact
that the boundaries of $G^{\varepsilon}$ have the same
regularity as that of $G$ if $\varepsilon$ is small enough.

Take a nonnegative symmetric $\zeta\in C^{\infty}_{0}
(B_{1})$ with unit integral, for $\varepsilon>0$
define $\zeta_{\varepsilon}(x)=\varepsilon^{-d}
\zeta(x/\varepsilon)$, and for functions $v$
on $\bR^{d}$ set $v^{(\varepsilon)}=\zeta_{\varepsilon}
*v$.  

\begin{lemma}
                                   \label{lemma 1.9.2}
Let $\varepsilon_{0}>0$ be sufficiently small, $p\in[1,
\infty)$, and 
$u\in W^{2}_{\infty}(G^{\varepsilon_{0}})$. 
Assume that $F[u]=0$ in $G^{\varepsilon_{0}}$.
Then for
any $h,\varepsilon\in(0,\varepsilon_{0}/4]$ 
there exists an $x_{0}$ with $|x_{0}|
\leq h$ for which
$$
\sum_{x\in   G_{h}}\big|F(u^{(\varepsilon)}(x+x_{0}),
Du^{(\varepsilon)}(x+x_{0}),D^{2}
u^{(\varepsilon)}(x+x_{0}),x)\big|^{p}h^{d}
$$
\begin{equation}
                                        \label{1.9.3}
\leq N(h^{p\gamma_{1}}+\varepsilon^{\gamma_{1}\gamma_{0}})
(1+M_{\varepsilon_{0}}^{p}),
\end{equation}
where $N$ depends only on $G$, $\delta$,   $p$,
and $d$.

\end{lemma}

Proof. The terms in the sum in \eqref{1.9.3} are obviously less
than $h^{d}$ times
$
2^{p}I(x+x_{0})+2^{p}J(x,x_{0})
$,
where
$$
I(x)=|F[u^{(\varepsilon)}](x)-F[u](x)|^{p},
$$
$$
J(x,x_{0})=
|F(u^{(\varepsilon)}(x+x_{0}),
Du^{(\varepsilon)}(x+x_{0}),D^{2}u^{(\varepsilon)}(x+x_{0}),x)
$$
$$
-F(u^{(\varepsilon)}(x+x_{0}),
Du^{(\varepsilon)}(x+x_{0}),D^{2}u^{(\varepsilon)}(x+x_{0}),x
+x_{0})|^{p}.
$$
By assumption
$
J(x,x_{0})\leq Nh^{p\gamma_{1}}(1+M_{\varepsilon_{0}}^{p})
$, so that
$$
\sum_{x\in G_{h}}J(x,x_{0})h^{d}
\leq Nh^{p\gamma_{1}}(1+M_{\varepsilon_{0}}^{p})
$$
whenever $|x_{0}|\leq h$.

Next, observe that
$$
\sum_{x\in  G_{h}}I_{x+B_{h/2}}\leq 
I_{G^{\varepsilon_{0}/2}}
$$
implying that
$$
\dashint_{ B_{h/2}}\sum_{x\in  G_{h}}
I(x+x_{0})h^{d}
\,dx_{0}=N\sum_{x\in  G_{h}}\int_{ B_{h/2}} 
I(x+x_{0})
\,dx_{0}
$$
$$
=\sum_{x\in  G_{h}}\int_{x+ B_{h/2}} 
I(y)
\,dy\leq\int_{G^{\varepsilon_{0}/2}}I(x)
\,dx.
$$

Furthermore, for $x\in G^{\varepsilon_{0}/2}$
we have
$$
|D^{2}u^{(\varepsilon)}(x)-D^{2}u(x)|^{p}\leq
N\dashint_{B_{\varepsilon}}|D^{2}u(x+y)-D^{2}u(x)|^{p}\,dy.
$$
Similar relations are true for the first
order derivatives and functions themselves.
Therefore, in light of Theorem  \ref{theorem 1.9.1}
$$
\int_{G^{\varepsilon_{0}/2}}I(x)\,dx
\leq N\dashint_{B_{\varepsilon}}
\|
u(\cdot+y)-u\|_{W^{2}_{p}(G^{\varepsilon_{0}/2})}^{p}
\,dy\leq
N_{1}\varepsilon^{\gamma_{1}\gamma_{0}}
(1+M^{p}_{\varepsilon_{0}}).
$$
We conclude  that
there exists an  $x_{0}$ with $|x_{0}|< h/2$ such that
$$
\sum_{x\in G_{h}}
I(x+x_{0})h^{d}\leq 2
N_{1}\varepsilon^{\gamma_{1}\gamma_{0}}
(1+M^{p}_{\varepsilon_{0}}).
$$

This along with the above estimate of $J$
proves the lemma. 

Below in the paper $\Lambda$ is not
necessarily the set from Section
\ref{section 1.30.1} and we introduce
$G^{o}$ and $\partial_{h}G$ in the same way
as in Section
\ref{section 1.30.1} for any $\Lambda$ we choose.

\begin{definition}
                           \label{definition 1.13.1}
Let $\Lambda\subset \bZ^{d} $.
We say that an operator
$H_{h}$ defined on functions on $\bZ^{d}_{h}$
and mapping them into functions on $\bZ^{d}_{h}$
is $\Lambda$-local if for any $x_{0}\in 
\bZ^{d}_{h}$ and  functions $u$ and $v$ 
on $\bZ^{d}_{h}$ such that $u=v$ on $x_{0}+h\Lambda$
we have $H_{h}[u](x_{0})=H_{h}[v](x_{0})$.
\end{definition}

\begin{definition}
                           \label{definition 1.12.1}
Let
$\cH=\{H_{h}:h\in(0,1)\}$ be a family
of $\Lambda$-local operators $H_{h}$ mapping functions on $\bZ^{d}_{h}$ to functions
on $\bZ^{d}_{h}$ and $p\geq1$. We say
that the family is $\Lambda\cL_{p}$-stable in $G$ if there 
exists constants $N=N(\cH,p)$ and $h_{0}=h_{0}(\cH,p)>0$
 such that, for any 
$h\in(0,h_{0}]$ and functions
$v $ and $u $ on $\bZ^{d}_{h}$, in $G_{h}$ we have
$$
u -v \leq N\bigg(\sum_{x\in  
 G^{o}_{h}}
\big(H_{h}[u](x)-H_{h}[v](x)\big)^{p}_{-}
h^{d}\bigg)^{1/p}+
\max_{\partial_{h}G}(u -v)_{+}.
$$
\end{definition}

\begin{lemma}
                                 \label{lemma 1.12.1}
Assume that $\cH$ is $\Lambda\cL_{p}$-stable in $G$. Also
assume that  
\begin{equation}
                                       \label{1.14.1}
F[u]-H_{h}[u]\geq -N_{1}h^{\gamma_{1}}\|u\|_{C^{3}(G)},
\end{equation}
in $G^{o}_{h}$ for any $u\in C^{3}(\bar{G})$ 
and $h\in(0,1)$
 or
\begin{equation}
                                       \label{1.28.9}
F[u]-H_{h}[u]\leq N_{1}h^{\gamma_{1}}\|u\|_{C^{3}(G)},
\end{equation}
in $G^{o}_{h}$ for any $u\in C^{3}(\bar{G})$ 
and $h\in(0,h_{0}]$,
where $N_{1}$ is a constant independent of $u$
and $h$. Finally, suppose that for each $h\in(0,h_{0}]$ we are
given a function $v_{h}$ on $G_{h}$ such that
$H_{h}[v_{h}]=0$ on $G^{o}_{h}$. Then for any 
sufficiently small $
\varepsilon_{0}>0 $ and $h\leq h_{0}$ satisfying
  $h^{\gamma_{1}/2}\leq\varepsilon_{0}/4$  and
$u\in W^{2}_{\infty}(G^{\varepsilon_{0}})$ such that
$F[u]=0$ in $G^{\varepsilon_{0}}$ we have 
that on $G_{h}$
\begin{equation}
                                      \label{1.14.2}
v_{h}\leq u+\max_{\partial_{h}G}(v_{h}
 -u)_{+}+Nh^{\alpha}(1+M_{\varepsilon_{0}})
\end{equation}
if \eqref{1.14.1} holds and
\begin{equation}
                                      \label{1.28.01}
u\leq v_{h}+\max_{\partial_{h}G}(u-v_{h}
 )_{+}+Nh^{\alpha}(1+M_{\varepsilon_{0}})
\end{equation}
if \eqref{1.28.9} holds,
where $N$ and $\alpha\in(0,\gamma_{1}/2)$ depend only on $N_{1}$,
$N(\cH,p)$,
$G$, $\delta$,  $p$,
and $d$. 
\end{lemma}

Proof. First assume that \eqref{1.14.1} holds.
By Lemma \ref{lemma 1.9.2}
for
any $h,\varepsilon\in(0,\varepsilon_{0}/4]$ 
there exists an $x_{0}$ with $|x_{0}|
\leq h$ for which \eqref{1.9.3} holds.
Owing to \eqref{1.14.1} and the fact that
$H_{h}[v_{h}]=0$,
if $x\in G_{h}^{o}$, we infer that
$$
\sum_{x\in  G^{o}_{h}}[H_{h}(u^{(\varepsilon)}
(\cdot+x_{0}))(x)-
H_{h}[v_{h}](x)]_{+}^{p}h^{d}
$$
$$
\leq N(h^{p\gamma_{1}}+\varepsilon^{\gamma_{1}\gamma_{0}})
(1+M_{\varepsilon_{0}}^{p})+NN^{p}_{1}h^{p\gamma_{1}}
\|u^{(\varepsilon)}\|^{p}_{C^{3}(G^{\varepsilon_{0}/4})}.
$$

According to Definition \ref{definition 1.12.1}, if $h\leq h_{0}$,
in $G_{h}$
$$
v_{h}\leq u^{(\varepsilon)}(\cdot+x_{0})
+ N(h^{ \gamma_{1}}+\varepsilon^{\gamma_{1}\gamma_{0}/p})
(1+M_{\varepsilon_{0}})+NN_{1}h^{\gamma_{1}}
\|u^{(\varepsilon)}\|_{C^{3}(G^{\varepsilon_{0}/4})}
$$
$$
+\max_{\partial_{h}G}(v_{h}
 -u^{(\varepsilon)}(\cdot+x_{0}))_{+}.
$$

One knows that
$$
\|u^{(\varepsilon)}\|_{C^{3}(G^{\varepsilon_{0}/4})}
\leq N\varepsilon^{-1}M_{\varepsilon_{0}}.
$$
Furthermore, for $x\in G$
$$
u^{(\varepsilon)}(x+x_{0})\geq
u(x+x_{0})-NM_{\varepsilon_{0}} \varepsilon^{2}\geq
u(x)-NM_{\varepsilon_{0}}(h+\varepsilon^{2})
$$
and $u^{(\varepsilon)}(x+x_{0})$ admits a similar estimate from
above.

Hence, in $G_{h}$ for
any $h,\varepsilon\in(0,\varepsilon_{0}/4]$ satisfying  
$h\leq h_{0}$ we have
$$
v_{h}\leq u
+\max_{\partial_{h}G}(v_{h}
 -u)_{+}
$$
$$
+NM_{\varepsilon_{0}}(h+\varepsilon^{2})
+N(h^{ \gamma_{1}}+\varepsilon^{\gamma_{1}\gamma_{0}/p})
(1+M_{\varepsilon_{0}})+NN_{1}h^{\gamma_{1}}
\varepsilon^{-1}M_{\varepsilon_{0}}
$$
and the result follows if we take $\varepsilon
=h^{\gamma_{1}/2}$ and restrict $h$ to satisfy
$h^{\gamma_{1}/2}\leq\varepsilon_{0}/4$
in order to have $\varepsilon\leq
\varepsilon_{0}/4$. 

This proves our assertion concerning 
\eqref{1.14.2}. The remaining assertion is
proved similarly.
The lemma is proved.

\mysection{Application of the general scheme}

Let $H(u,x)$ be a function satisfying Assumption
\ref{assumption 1.14.1} and such that
\begin{equation}
                                                \label{1.21.2}
|H(0,x)|\leq\delta^{-1}.
\end{equation}
We start with the following consequences of 
Theorems 1.1 and 1.3 of \cite{Kr_12}.

\begin{theorem}
                             \label{theorem 9.23.1}
There is a  constant  $\hat{\delta}\in(0,\delta]$
depending only on $\delta$ and $d$ and there exists
a function $P(u) $ (independent of $x$), 
satisfying Assumption
\ref{assumption 1.14.1} 
with $\hat{\delta}$   
in  place of 
$\delta$, 
such that for any $K\geq0$ and sufficiently small
$\varepsilon_{0}\geq0$  each of the equations
\begin{equation}
                                       \label{9.23.2}
\max(H[v],P[v]-K)=0
\end{equation}
and
\begin{equation}
                                       \label{1.28.4}
\min(H[v],-P[-v]+K)=0
\end{equation}
in $G^{\varepsilon_{0}}$  (a.e.)
 with boundary condition $v=0$ on $\partial
G^{\varepsilon_{0}}$
has a unique solution in the space
$ C^{0,1}(\bar{G}^{\varepsilon_{0}})
\cap C^{1,1}_{loc}(G^{\varepsilon_{0}})$.
Denote by $v_{\varepsilon_{0},K}$ the solution of
\eqref{9.23.2} and by $v_{\varepsilon_{0},-K}$
 the solution of
\eqref{1.28.4}. Then 
\begin{equation}
                                                \label{1.13.1}
|v_{\varepsilon_{0},\pm K}|,
|Dv_{\varepsilon_{0},\pm K}|,\rho_{\varepsilon_{0}}
|D^{2}v_{\varepsilon_{0},\pm K} |\leq N(1+K)\quad
\text{in}\quad G^{\varepsilon_{0}}\quad (a.e.), 
\end{equation}
where 
$$
\rho_{\varepsilon_{0}}
(x)=\dist(x,\bR^{d}\setminus G^{\varepsilon_{0}}),
$$
and
   $N$ is a constant depending only on $G$ 
 and $\delta$.

Finally, $P(u )$ is
constructed on the sole basis of $\delta$ 
and $d$, it is  positive homogeneous of degree one
and convex in $u$.

\end{theorem}

Actually, Theorems 1.1 and 1.3 of \cite{Kr_12} are proved
only for $\varepsilon_{0}=0$. The fact that
they also hold  for sufficiently small
$\varepsilon_{0}>0$ easily follows by inspecting
their proofs.

\begin{theorem}
                              \label{theorem 1.18.1}
Assume that we are given a  $\Lambda\cL_{p}$-stable
family $\cH$ of operators $H_{h}$. Also suppose that
\begin{equation}
                                  \label{1.28.7}
H[u]-H_{h}[u]\geq -N_{1}h^{\gamma_{1}}\|u\|_{C^{3}(G)},
\end{equation}
in $G^{o}_{h}$ for any $u\in C^{3}(\bar{G})$ 
and $h\in(0,h_{0}]$ or
\begin{equation}
                                  \label{1.28.6}
H[u]-H_{h}[u]\leq N_{1}h^{\gamma_{1}}\|u\|_{C^{3}(G)},
\end{equation}
in $G^{o}_{h}$ for any $u\in C^{3}(\bar{G})$ 
and $h\in(0,h_{0}]$,
where $N_{1}$ is a constant independent of $u$
and $h$. Finally, suppose that for each $h\in(0,h_{0}]$ we are
given a function $v_{h}$ on $G_{h}$ such that
$H_{h}[v_{h}]=0$ on $G^{o}_{h}$ and 
$v_{h}=0$ on $\partial_{h}G$. Then 
there exists an $h_{1}\in(0,h_{0}]$ depending only
on $G$, $h_{0}$, $\alpha$ and $\delta$ such that for
any $K\geq0$ and $h\in(0, h_{1}]$ on $G_{h}$
we have
\begin{equation}
                                  \label{1.19.01}
v_{h}\leq v_{0,K} 
+N(1+K)h^{\alpha/2},
\end{equation}
if \eqref{1.28.7} holds and 
\begin{equation}
                                  \label{1.28.8}
v_{h}\geq v_{0,-K} 
+N(1+K)h^{\alpha/2}
\end{equation}
in case  \eqref{1.28.6} holds,
where $N$ depends only on
$\Lambda$, $N_{1}$,
$N(\cH,p)$,
$G$, $\delta$,  $p$,
and $d$ and $\alpha$
is the constant from Lemma \ref{lemma 1.12.1}.

\end{theorem}

Proof. First suppose that 
\eqref{1.28.7} holds. For $K$ fixed introduce $F(u,x)
=\max(H(u,x),P(u)-K)$ and observe that
$F\geq H$, so that \eqref{1.28.7}
implies \eqref{1.14.1}. Therefore,
we can apply Lemma \ref{lemma 1.12.1}
with $v_{\varepsilon_{0},K}$ for
$\varepsilon_{0}>0$ in place of $u$.

Notice that $\rho_{\varepsilon_{0}}\geq
\varepsilon_{0}$ on $G$, so that
$M_{\varepsilon_{0}}\leq N(1+K)
\varepsilon_{0}^{-1}$ by Theorem
\ref{theorem 9.23.1}. By the same theorem
$| v_{\varepsilon_{0},K}|\leq N\varepsilon_{0}
(1+K)$ in $\partial_{h}G$ if $h 
\leq\varepsilon_{0}$. Hence, for  $h$ also satisfying
$h\leq h_{0}$ and $4h^{\gamma_{1}/2}\leq\varepsilon_{0}$
and sufficiently small
$\varepsilon_{0}$,  on $G_{h}$
\begin{equation}
                                   \label{1.19.1}
v_{h}\leq v_{\varepsilon_{0},K}+N\varepsilon_{0}(1+K)
+Nh^{\alpha}(1+\varepsilon_{0}^{-1}(1+K)).
\end{equation}

Next, for a constant $N$ depending only on $G$
and $\delta$ we have that $v_{\varepsilon_{0},K}
\leq N\varepsilon_{0}
(1+K)=v_{0,K}+ N\varepsilon_{0}
(1+K)$ on $\partial G$. By the maximum principle
this inequality extends to $G$ and in light of
 \eqref{1.19.1} yields that  on $G_{h}$
\begin{equation}
                                   \label{1.19.3} 
v_{h}\leq v_{0,K}+N\varepsilon_{0}(1+K)
+Nh^{\alpha}(1+\varepsilon_{0}^{-1}(1+K)).
\end{equation}
We take here $\varepsilon_{0}=4h^{\alpha/2}$ and call
$h_{1}$ the largest value of $h\in(0,h_{0}]$ 
for which $\varepsilon_{0}=4h_{1}^{\alpha/2}$
is sufficiently small and less than 1. Observe that
then $4h^{\gamma_{1}/2}\leq\varepsilon_{0}$ for 
$h\leq h_{1}$
because $\alpha\leq\gamma_{1}$ and $h_{1}\leq1$. Now
  \eqref{1.19.3} yields \eqref{1.19.01}.

  The assertion
about \eqref{1.28.8} is proved similarly.
The theorem is proved.

\begin{corollary}
                                      \label{corollary 1.21.1}
Assume that we are given a  $\Lambda\cL_{p}$-stable
family $\cH$ of operators $H_{h}$. Also suppose that
\begin{equation}
                                                      \label{1.30.1}
|H[u]-H_{h}[u]|\leq N_{1}h^{\gamma_{1}}\|u\|_{C^{3}(G)},
\end{equation}
in $G^{o}_{h}$ for any $u\in C^{3}(\bar{G})$ 
and $h\in(0,h_{0}]$,
where $N_{1}$ is a constant independent of $u$
and $h$. Suppose that for each $h\in(0,h_{0}]$ we are
given a function $v_{h}$ on $G_{h}$ such that
$H_{h}[v_{h}]=0$ on $G^{o}_{h}$ and 
$v_{h}=0$ on $\partial_{h}G$.
Finally, suppose that there exist constants $N_{2}$ and $\gamma_{2}>0$
such that  for the unique viscosity
solution $v\in C(\bar{G})$
  of equation $H[v]=0$ in $G$ with zero boundary data
we have 
\begin{equation}
                                                      \label{1.28.1}
v_{0, - K}+N_{2}K^{-\gamma_{2}}\geq 
v\geq v_{0,  K}-N_{2}K^{-\gamma_{2}}
\end{equation}
in $G$ for any $K\geq1$. Then 
there exist constants   $\beta\in(0,1)$ and $N$ such that for
all sufficiently small $h>0$ we have that on $G_{h}$
$$
|v_{h}-v|\leq Nh^{\beta}.
$$

\end{corollary}

Indeed, we have for $K\geq1$ and $h\in(0,h_{1}]$ that
 on $G_{h}$
$$
v_{h}\leq v_{0,K} 
+N(1+K)h^{\alpha/2}\leq v  
+N(1+K)h^{\alpha/2}+N_{2}K^{-\gamma_{2}}
$$
and to show that 
$v_{h}\leq v+ Nh^{\beta}$,
it only remains to minimize the last expression
with respect to $K\geq1$. The estimate
from the other side is obtained similarly.

\mysection{Proof of Theorem \protect\ref{theorem 1.28.1}}
                                   \label{section 1.30.2}

Obviously, $H$ is Lipschitz continuous with respect to $u$
with Lipschitz constant depending only on $\delta$ and $d$.
This function is also decreasing with respect to $u'_{0}$
and therefore $F_{u'_{0}} \leq 0$ wherever $F_{u'_{0}}$
exists. Furthermore, for $t>0$ and $x,\lambda\in\bR^{d}$
we have
$$
H(u',u''+t\lambda\lambda^{*},x)\geq H(u,x)+t\delta|\lambda|^{2},
$$
which shows that $H_{u''}$, whenever it exists,
is uniformly nondegenerate, so that 
Assumption \ref{assumption 1.14.1} (ii)
is satisfied with a perhaps different $\delta>0$.
Assumption \ref{assumption 1.14.1} (iii)
is satisfied as well with a perhaps  different $\delta>0$.

The properties of $P$ listed in Theorem \ref{theorem 9.23.1}
or just the construction of $P$ in \cite{Kr_12}
yield that there is a set $A_{2}$ and bounded 
continuous
functions $\sigma^{\alpha}=\sigma^{\alpha\beta}$,
$b^{\alpha}=b^{\alpha\beta}$, $c^{\alpha}=c^{\alpha\beta}$  
(independent of  $x$ and $\beta$), 
and $f^{\alpha\beta}\equiv0$ defined on $A_{2}$ 
such that   Assumption \ref{assumption 1.9.1}
is satisfied perhaps with a different constant $\delta>0$
and for $a^{\alpha}:=a^{\alpha\beta}
=(1/2)\sigma^{\alpha}(\sigma^{\alpha})^{*}$ we have
\begin{equation}
                                                  \label{2.8.5}
P[u](x)=\sup_{\alpha\in A_{2}}
\big[a_{ij}^{\alpha} D_{ij}u(x)
+b^{\alpha}_{i} D_{i}u(x)-c^{\alpha} u(x) \big].
\end{equation}
Define $A_{1}=A$,
$$
\hat{A}=A_{1}\cup A_{2}
$$
and observe that
$$
\max(H[u](x),
P[u](x)-K)
$$
$$
=\max\big\{\supinf_{\alpha\in A_{1}\,\,\beta\in B}
[L^{\alpha\beta}u(x)+f^{\alpha\beta}(x)],
\supinf_{\alpha\in A_{2}\,\,\beta\in B}
[L^{\alpha\beta}u(x)+f^{\alpha\beta}(x)-K]\big\}
$$
$$
=\supinf_{\alpha\in\hat{A}\,\,\beta\in B}
\big[ L^{\alpha\beta}u( x)+f^{\alpha\beta}_{K}( x)]
\quad (f^{\alpha\beta}_{K}( x)=f^{\alpha\beta} ( x)
I_{\alpha\in A_{1}}-KI_{\alpha\in A_{2}}),
$$
where the first equality follows from the definition of $H[u]$,
\eqref{2.8.5}, and the fact that $L^{\alpha \beta}$
is independent of $\beta$ for $\alpha\in A_{2}$.

We have just repeated part of the proof of Theorem 2.3
of \cite{Kr_14}. Then, as there, we have a probabilistic representation
of $v_{0,K}$, which
allows us to conclude that
$|v_{0,K}-v|\leq N_{2}K^{-\gamma_{2}}$ in $G$
by inspecting the proof
of Theorem 5.2 of \cite{Kr_14} in which the convergence
$v_{0,K}\to v$ is proved under much weaker assumptions
on $c^{\alpha\beta}$ and $f^{\alpha\beta}$.
Exploiting the H\"older continuity of
$c^{\alpha\beta}$ and $f^{\alpha\beta}$ easily
yields the stated rate of convergence of 
$v_{0,K}$ to $v$.

 Passing to $v_{0,-K}$ denote $B_{1}=B$, $B_{2}=A_{2}$,
$$
\hat B=B_{1}\cup B_{2}
$$
and notice that 
$$
-P[-u](x)=\inf_{\beta\in B_{2}}
\big[a_{ij}^{\beta}D_{ij}u(x)
+b^{\beta}_{i}D_{i}u(x)-c^{\beta}u(x)\big].
$$
Next, as in \cite{Kr_14}
$$
 \min(H[u](x),
-P[-u](x)+K) 
$$
$$
=\sup_{\alpha\in A }\min\big\{\inf_{ \beta\in B_{1}}
[  L^{\alpha\beta}u(x)+ f^{\alpha\beta}(x)],
 \inf_{\beta\in B_{2}}
[L^{\alpha\beta}u(x)+f^{\alpha\beta}(x)+K]\big\}
$$
$$
 =\supinf_{\alpha\in A\,\,\beta\in \hat{B}}
\big[ L^{\alpha\beta}u( x)+f^{\alpha\beta}_{K}( x)]
$$
with $f^{\alpha\beta}_{K}$ defined this time by
$$
f^{\alpha\beta}_{K}( x)=f^{\alpha\beta}( x)
I_{\beta\in B_{1}}+KI_{\beta\in B_{2}}
$$
and $L^{\alpha\beta}u=
a_{ij}^{\beta}D_{ij}u 
+b^{\beta}_{i}D_{i}u -c^{\beta}u $ for $\beta\in B_{2}$.
This allows us to get a probabilistic representation
of $v_{0,-K}$,  which
allows us to conclude that
$|v_{0,-K}-v|\leq N_{2}K^{-\gamma_{2}}$ in $G$
by inspecting this time the proof
of Theorem 6.1 of \cite{Kr_14}.
We thus checked the assumption \eqref{1.28.1}.

 Furthermore,
it follows from Taylor's formula that in $G^{0}_{h}$ we have
$$
|L^{\alpha\beta}u-L^{\alpha\beta}_{h}u|\leq Nh\|u\|_{C^{3}(G)}
$$
for any $u\in C^{3}(\bar{G})$, where $N$ is independent of
$\alpha,\beta$ and $h$. We see that condition
\eqref{1.30.1} is satisfied.

Finally, as is shown in  \cite{KT92}, the family of operators
$H_{h}$ is $\Lambda\cL_{d}$-stable and now
a simple reference to
 Corollary \ref{corollary 1.21.1}
proves our  Theorem \ref{theorem 1.28.1}.

\end{document}